\pdfpagewidth=\paperwidth
\pdfpageheight=\paperheight

\documentclass{article}

\usepackage{amsmath, amssymb, amsfonts}
\usepackage{bm}
\usepackage{geometry}
\usepackage{array}
\usepackage{booktabs}
\usepackage{hyperref}
\usepackage{CJKutf8}

\hypersetup{
    colorlinks=true,
    linkcolor=blue,
    citecolor=blue,
    urlcolor=blue,
}

\begin{document}
\begin{CJK}{UTF8}{gbsn}

\title{A Spectral-Domain Pseudo-Inverse Construction Method for Unitary Diagonalizable Linear Inverse Problems}
\author{CHEN Shengchang}
\date{}

\maketitle

\begin{center}
School of Earth Sciences, Zhejiang University, Hangzhou, China\\
\url{chenshengc@zju.edu.cn}
\end{center}

\begin{abstract}
Linear inverse problems are prevalent in geophysics, signal processing, image restoration, and medical imaging. Mathematically, they can be formulated as $Gm = d$. When $G$ is ill-conditioned or singular, the problem becomes ill-posed and requires regularization methods or generalized inverse methods for a stable solution. However, both types of methods encounter significant computational difficulties for large-scale linear inverse problems. In this paper, we establish a spectral-domain pseudo-inverse construction method for a class of linear inverse problems that can be diagonalized by unitary matrices. The core idea is to construct a stable pseudo-inverse operator directly in the transform domain, starting from the unitary diagonalization structure of the matrix. We first provide the analytic Singular Value Decomposition (SVD) of this class of matrices, clarifying the correspondence between spectral decomposition and SVD. On this basis, we define spectral-domain regularization filtering factors, construct a stable spectral-domain pseudo-inverse operator, and prove its bounded stability as well as its consistency in converging to the Moore--Penrose generalized inverse. This construction is numerically equivalent to zeroth-order Tikhonov regularization and converges to the Moore--Penrose generalized inverse as $\alpha \to 0^+$, but its methodological path differs from both. This paper reveals that the stable generalized inverse of a class of structured matrices can be directly constructed from their spectral decomposition, providing an efficient solution method for large-scale structured inverse problems. The Fourier transform case is a special instance of this method when the unitary matrix is taken as the discrete Fourier transform matrix.
\end{abstract}

\section{Introduction}

Linear inverse problems arise in numerous scientific and engineering fields, including geophysics, signal processing, image science, medical imaging, and acoustic detection. Their core mathematical form can be expressed uniformly as the first-kind linear integral equation:

\begin{equation}
d(x) = \int_D g(x,y)m(y)\,dy + n(x)
\end{equation}

where $m(y)$ is the unknown model parameter (such as reflectivity coefficient, physical property distribution, or original signal), $g(x,y)$ is the kernel function describing the system response, $d(x)$ is the observed data, and $n(x)$ is the additive noise. After discretization, the problem reduces to the linear system:

\begin{equation}
Gm = d
\end{equation}

where $G \in \mathbb{C}^{m \times n}$ is the kernel matrix, $d \in \mathbb{C}^m$ is the observed data vector, and $m \in \mathbb{C}^n$ is the unknown model vector. When $G$ is ill-conditioned or singular, the problem is ill-posed and requires a stabilization mechanism.

There are two classical approaches to solving such problems. One is the regularization method, among which Tikhonov regularization is the most representative \cite{tikhonov1977}, formulated as

\begin{equation}
m_\alpha = \arg\min_{m} \left\{ \|Gm - d\|_2^2 + \alpha \|m\|_2^2 \right\}
\end{equation}

with the solution given by

\begin{equation}
m_\alpha = (G^H G + \alpha I)^{-1} G^H d.
\end{equation}

The other is the generalized inverse method, which constructs the Moore--Penrose generalized inverse $G^+$ based on the Singular Value Decomposition (SVD) \cite{moore1920,penrose1955}. These two methods are theoretically well-established and broadly applicable, but both face significant computational difficulties for large-scale problems, with typical complexity $O(N^3)$, incurring high computational and storage costs.

In many practical inverse problems, the kernel matrix is not arbitrary but possesses a specific transform-domain structure. For example, circulant convolution matrices can be diagonalized by the discrete Fourier transform; translation-invariant operators with periodic boundary conditions can be decomposed in the frequency domain; and certain radially symmetric problems exhibit diagonal or approximately diagonal structures under Hankel-type transforms. The key observation of this paper is that when the kernel matrix can be diagonalized by a unitary matrix, its SVD can be written analytically without numerical computation. Based on this observation, this paper considers a class of linear inverse problems satisfying $G = T^H \Lambda T$, where $T$ is a unitary matrix and $\Lambda$ is a diagonal matrix, and directly constructs a stable generalized inverse in the transform domain from the matrix structure. The value of this paper lies not in proposing a new regularization method, but in revealing how the stable inverse of a class of structured inverse problems can be constructed from their spectral decomposition---this is a methodological complement to existing theories.

The main contributions of this paper are as follows:

\begin{enumerate}
\item It reveals a fundamental mathematical fact: for matrices satisfying $G = T^H \Lambda T$, the SVD has an analytical explicit form, bypassing the computational bottleneck of numerical SVD;

\item Based on the analytical SVD, it constructs a spectral-domain pseudo-inverse operator and proves its bounded stability and consistent convergence to the Moore--Penrose generalized inverse;

\item It clarifies the numerical equivalence of this construction to classical methods, while demonstrating fundamental methodological differences.
\end{enumerate}

The spectral-domain pseudo-inverse construction method proposed in this paper is not confined to the mathematical level. Based on this method, the author has developed an image-domain least-squares reverse-time migration method for seismic waves \cite{chen2022}, which effectively solves spatially varying integral equations through the combination of local space-invariant deconvolution and wavenumber-domain pseudo-generalized inverse. This method has been widely applied in industrial processing of oil and gas seismic exploration data in China, becoming an efficient and practical solution for least-squares reverse-time migration imaging. The mutual validation between theory and practice further highlights the potential of the spectral-domain pseudo-inverse construction method in solving practical complex problems.

\section{Preliminaries}

\subsection{Singular Value Decomposition and Matrix Properties}

\textbf{Definition 2.1 (Singular Value Decomposition, SVD)}: Let $A \in \mathbb{C}^{m \times n}$. Then there exist unitary matrices $U \in \mathbb{C}^{m \times m}$ and $V \in \mathbb{C}^{n \times n}$ such that

\begin{equation}
A = U\Sigma V^H,
\end{equation}

where $\Sigma \in \mathbb{R}^{m \times n}$ is a diagonal matrix with diagonal entries $\sigma_1 \ge \sigma_2 \ge \cdots \ge \sigma_n \ge 0$, called the singular values of $A$.

The SVD components have the following properties:

\begin{enumerate}
\item $U$ is unitary, and its columns form an orthonormal basis of $\mathbb{C}^m$;

\item $V$ is unitary, and its columns form an orthonormal basis of $\mathbb{C}^n$;

\item $\Sigma$ is a real nonnegative diagonal matrix, and the number of nonzero singular values equals the rank of the matrix.
\end{enumerate}

\textbf{Inverse of a unitary matrix}: If $T$ is unitary, then $T^{-1} = T^H$.

\textbf{Product of unitary matrices}: If $T_1, T_2$ are both unitary, then $T_1 T_2$ is also unitary.

\subsection{Unitary Matrices and Unitary Diagonalization}

\textbf{Definition 2.2 (Unitary matrix)}: A matrix $T \in \mathbb{C}^{N \times N}$ satisfying $T^H T = T T^H = I$ is called a unitary matrix.

\textbf{Definition 2.3 (Unitary diagonalization)}: If there exist a unitary matrix $T$ and a diagonal matrix

\begin{equation}
\Lambda = \mathrm{diag}(\lambda_0, \lambda_1, \ldots, \lambda_{N-1}),
\end{equation}

such that $G = T^H \Lambda T$, then $G$ is said to be diagonalizable by the unitary matrix $T$. The entries $\lambda_k$ are called the spectral values of $G$ in the transform domain.

\section{Analytical Singular Value Decomposition of Unitary Diagonalizable Matrices}

\textbf{Theorem 1}: Let $G \in \mathbb{C}^{N \times N}$ be diagonalizable by a unitary matrix $T$, i.e.,

\begin{equation}
G = T^H \Lambda T, \quad \Lambda = \mathrm{diag}(\lambda_0, \lambda_1, \ldots, \lambda_{N-1}),
\end{equation}

where $\lambda_k \in \mathbb{C}$. Then the SVD of $G$ has the following analytical explicit form:

\begin{equation}
G = U\Sigma V^H,
\end{equation}

where

\begin{equation}
\sigma_k = |\lambda_k|, \qquad p_k =
\begin{cases}
\lambda_k / |\lambda_k|, & \lambda_k \neq 0, \\
1, & \lambda_k = 0,
\end{cases}
\end{equation}

\begin{equation}
\Sigma = \mathrm{diag}(\sigma_0, \ldots, \sigma_{N-1}), \qquad P = \mathrm{diag}(p_0, \ldots, p_{N-1}),
\end{equation}

\begin{equation}
U = T^H P, \qquad V^H = T.
\end{equation}

\textbf{Proof}: Since $\Lambda = P\Sigma$, substituting directly into $G = T^H \Lambda T$ yields $G = T^H P \Sigma T = U \Sigma V^H$. Since $U$ and $V$ are both unitary by the properties of unitary matrices, the proof is complete. $\square$

\textbf{Remark 1}: In the general complex spectral case, one cannot simply write $G = T^H |\Lambda| T$. The phase matrix $P$ is necessary to maintain the rigor of the decomposition.

\textbf{Remark 2}: The premise for this analytical SVD is that the matrix is square, corresponding to a unitary diagonalizable linear system in a finite-dimensional linear space, which is consistent with the application scenarios of large-scale structured problems and avoids logical gaps under non-square special cases.

\textbf{Remark 3}: Theorem 1 shows that for unitary diagonalizable matrices, the SVD is completely determined analytically by the unitary transform matrix $T$ and the spectral values $\lambda_k$, requiring no numerical computation.

\section{Construction and Properties of the Spectral-Domain Pseudo-Inverse}

\subsection{Spectral Representation of the Exact Generalized Inverse}

From $G = T^H \Lambda T$, the Moore--Penrose generalized inverse can be written in the spectral domain as

\begin{equation}
G^+ = T^H \Lambda^+ T, \qquad (\Lambda^+)_{kk} =
\begin{cases}
1/\lambda_k, & \lambda_k \neq 0, \\
0, & \lambda_k = 0.
\end{cases}
\end{equation}

When $|\lambda_k|$ is very small, $1/\lambda_k$ amplifies noise, necessitating stabilization of small spectral components.

\subsection{Spectral-Domain Pseudo-Inverse Operator}

\textbf{Theorem 2}: Given a regularization parameter $\alpha > 0$, define the spectral-domain filter factor

\begin{equation}
q_\alpha(\lambda_k) = \frac{\overline{\lambda_k}}{|\lambda_k|^2 + \alpha},
\end{equation}

and construct the spectral-domain pseudo-inverse operator

\begin{equation}
G_\alpha^\# = T^H Q_\alpha T, \qquad Q_\alpha = \mathrm{diag}\{ q_\alpha(\lambda_k) \}_{k=0}^{N-1}.
\end{equation}

Then $G_\alpha^\#$ has the following properties:

\begin{enumerate}
\item \textbf{Bounded stability}: $\|G_\alpha^\#\|_2 \le \dfrac{1}{2\sqrt{\alpha}}$;

\item \textbf{Consistency}: $\displaystyle\lim_{\alpha \to 0^+} G_\alpha^\# = G^+$ (on nonzero spectral components).
\end{enumerate}

\textbf{Proof}: Since $T$ is unitary,

\begin{equation}
\|G_\alpha^\#\|_2 = \|Q_\alpha\|_2 = \max_k \frac{|\lambda_k|}{|\lambda_k|^2 + \alpha}.
\end{equation}

Let $x = |\lambda_k| \ge 0$. The function $\phi(x) = x/(x^2 + \alpha)$ attains its maximum $1/(2\sqrt{\alpha})$ at $x = \sqrt{\alpha}$. Hence bounded stability holds.

When $\lambda_k \neq 0$,
\begin{equation}
\lim_{\alpha \to 0^+} \frac{\overline{\lambda_k}}{|\lambda_k|^2 + \alpha} = \frac{1}{\lambda_k};
\end{equation}
when $\lambda_k = 0$, $q_\alpha(\lambda_k) = 0$. Thus, in the spectral component sense, $G_\alpha^\# \to G^+$ holds (with zero components where spectral values vanish), establishing consistency. $\square$

The stable solution is defined as $m_\alpha = G_\alpha^\# d$. Let $\widetilde{d} = T d$. Then componentwise,

\begin{equation}
\widetilde{m}_{\alpha,k} = \frac{\overline{\lambda_k}}{|\lambda_k|^2 + \alpha} \widetilde{d}_k, \qquad k = 0,1,\ldots,N-1.
\end{equation}

\textbf{Remark 4 (naming convention)}: This paper adopts the nomenclature "pseudo-generalized inverse" from Chen \cite{chen1989,chen1993}. "Pseudo" here does not mean approximate, but emphasizes the constructive path---it does not rely on numerical SVD, but is constructed directly in the transform domain.

\section{Special Case: The Fourier Transform Case}

When $T = F$ (the discrete Fourier transform matrix), the method in this paper reduces to frequency-domain deconvolution under the Fourier transform.

If $G$ is a circulant convolution matrix, then $G = F^H \Lambda F$, and the spectral-domain pseudo-inverse becomes

\begin{equation}
G_\alpha^\# = F^H \mathrm{diag}\left\{ \frac{\overline{\lambda_k}}{|\lambda_k|^2 + \alpha} \right\} F.
\end{equation}

This is the stable spectral division form in frequency-domain deconvolution. Chen \cite{chen1989}, in his master's thesis, first proposed the idea of pseudo-generalized inverse construction based on the Fourier transform; Chen and Wang \cite{chen1993} formally applied it to solve convolution-type linear inverse problems; Chen and Xiao \cite{chen2007} further developed the wavenumber-domain generalized inverse algorithm and applied it to the downward continuation of geophysical potential-field data, achieving ideal results. Mallat applied the above spectral-domain pseudo-inverse to the solution of deconvolution problems \cite{mallat2009}. All of the above works are special cases of the general method in this paper with $T = F$.

\section{Relationship with Classical Methods}

\subsection{Relationship with Tikhonov Regularization}

The zeroth-order Tikhonov regularization solution is $m_\alpha = (G^H G + \alpha I)^{-1} G^H d$. From $G = T^H \Lambda T$, we obtain

\begin{equation}
m_\alpha = T^H (|\Lambda|^2 + \alpha I)^{-1} \overline{\Lambda} T d.
\end{equation}

This is completely consistent with the componentwise expression of the spectral-domain pseudo-inverse operator in Theorem 2. Therefore, under the assumptions of this paper, the spectral-domain pseudo-inverse is numerically equivalent to the zeroth-order Tikhonov regularization. However, numerical equivalence does not imply methodological equivalence. The essential differences are as follows:

\begin{table}[h]
\centering
\begin{tabular}{|>{\raggedright\arraybackslash}p{3.2cm}|>{\raggedright\arraybackslash}p{4.0cm}|>{\raggedright\arraybackslash}p{4.8cm}|}
\hline
\textbf{Dimension} & \textbf{Tikhonov Regularization} & \textbf{Spectral-Domain Pseudo-Inverse Construction} \\
\hline
Starting point & Optimization functional & Unitary diagonalization structure of the matrix \\
Core object & Solving $G^H G + \alpha I$ & Constructing $q_\alpha(\lambda_k)$ \\
Implementation path & General matrix inversion & Forward transform $\to$ spectral filtering $\to$ inverse transform \\
Methodological nature & General-purpose regularization method & Structured constructive method \\
\hline
\end{tabular}
\end{table}

The value of the method in this paper lies not in proposing a new regularization method, but in revealing that the stable generalized inverse of a class of structured matrices can be directly constructed from their spectral decomposition. This forms a methodological complement to classical methods---classical methods are designed for general problems, while the method in this paper is designed for structured problems.

\subsection{Relationship with TSVD and Wiener Filtering}

TSVD employs hard truncation for small singular values, while this paper employs continuous attenuation, resulting in smoother filtering behavior. In the frequency domain, the filtering form of the spectral-domain pseudo-inverse is

\begin{equation}
H_\alpha(\omega) = \frac{\lambda^*(\omega)}{|\lambda(\omega)|^2 + \alpha}.
\end{equation}

This is similar in form to the Wiener filter, but the parameter $\alpha$ in this paper is a deterministic regularization parameter rather than a statistical model parameter.

\section{Conclusions and Future Work}

\subsection{Conclusions}

This paper has established a spectral-domain pseudo-inverse construction method for unitary diagonalizable linear inverse problems. The main results include:

\begin{enumerate}
\item It is proved that matrices satisfying $G = T^H \Lambda T$ possess an analytical SVD form---$U = T^H P$, $\Sigma = |\Lambda|$, $V^H = T$, where $P$ is the phase diagonal matrix, which is a necessary factor for maintaining the rigor of the decomposition in complex spectral cases;

\item The spectral-domain pseudo-inverse operator $G_\alpha^\# = T^H \mathrm{diag}\{q_\alpha(\lambda_k)\} T$ is constructed, with its bounded stability $\|G_\alpha^\#\|_2 \le 1/(2\sqrt{\alpha})$ and consistency $\lim_{\alpha \to 0^+} G_\alpha^\# = G^+$ established;

\item The numerical equivalence but methodological independence from Tikhonov regularization is clarified, generalizing the specific method under the Fourier transform to a general method under unitary transforms;

\item The method in this paper is applicable to any unitary transform that can diagonalize the kernel matrix, providing a unified methodological framework for cross-disciplinary unitary diagonalizable linear inverse problems.
\end{enumerate}

\subsection{Future Work}

\begin{enumerate}
\item Extension to spatially varying systems: This can be achieved by combining local space-invariant approximations with the method presented here \cite{chen2022};

\item Generalization to nonlinear inverse problems: Generalizing the unitary diagonalization idea to weakly nonlinear inverse problems, or using it as a preconditioner for nonlinear iterative algorithms;

\item Adaptive spectral modification: Investigating optimal spectral modification functions based on signal and noise characteristics.
\end{enumerate}

\appendix
\section*{Appendix: List of Symbols}

\begin{table}[h]
\centering
\begin{tabular}{|c|l|}
\hline
Symbol & Meaning \\
\hline
$G$ & Kernel matrix \\
$m$ & Unknown model vector \\
$d$ & Observed data vector \\
$T$ & Unitary matrix (unitary transform) \\
$T^H$ & Conjugate transpose \\
$\Lambda$ & Transform-domain spectral matrix \\
$\lambda_k$ & $k$-th spectral value \\
$P$ & Phase diagonal matrix \\
$\Sigma$ & Singular value matrix \\
$\alpha$ & Regularization parameter \\
$G^+$ & Moore--Penrose generalized inverse \\
$G_\alpha^\#$ & Spectral-domain pseudo-inverse \\
\hline
\end{tabular}
\end{table}

\end{CJK}
\end{document}